\documentclass[12pt,a4paper,draft]{article}
\usepackage[T2A]{fontenc}
\usepackage{latexsym,amsfonts,amssymb,amsmath,longtable,amsthm,
mathrsfs}
\renewcommand{\le}{\leqslant}
\renewcommand{\ge}{\geqslant}

\newtheorem{Lemma}{{\bfseries Lemma}}
\newtheorem{Cor}[Lemma]{{\bfseries Corollary}}
\newtheorem{Theo}[Lemma]{{\bfseries Theorem}}
\newtheorem{Conj}[Lemma]{{\bfseries Conjecture}}

\DeclareMathOperator{\GL}{GL} 
\DeclareMathOperator{\SL}{SL}

\title{\vspace{-1cm} \hfill{\normalsize MSC 20D20}{
\fontfamily{cmr} \fontseries{bx} \selectfont \\ \vspace{1cm} On the inheriting
of the property $C_\pi$ by some normal subgroups}
\thanks{The work is supproted by RFBR, project 05-01-00797, grant of President
of RF for young scientists MK-3036.2007.1, and SB RAS Integration project
2006.1.2}}
\date{}
\author{\bf  E.P. Vdovin,
D.O. Revin}

\begin{document}



\maketitle
\pagenumbering{arabic}
\begin{abstract}

In the paper we prove that the  Hall property $C_\pi$ is inherited by normal
subgroups which index is a  $\pi'$-number.

\end{abstract}


\section*{Introduction}

Let  $\pi$ be a set of primes. We denote by  $\pi'$ the set of all primes not in
$\pi$, by $\pi(n)$ the set of prime divisors of a positive integer $n$, while
for a finite group  $G$ by $\pi(G)$ we denote $\pi(|G|)$. A positive integer
$n$ with $\pi(n)\subseteq\pi$ is called a $\pi$-number, whole a group $G$
$\pi(G)\subseteq \pi$ is called a $\pi$-group. A subgroup $H$ of $G$ is called
a {\it $\pi$-Hall subgroup}, if $\pi(H)\subseteq\pi$ by $\pi(|G:H|)\subseteq
\pi'$.  According to \cite{Hall} we shall say that $G$ {\it has the  property
$E_\pi$} (or briefly $G\in E_\pi$), if $G$ contains a $\pi$-Hall subgroup. If we
also have that every two $\pi$-Hall subgroups are conjugate, then we shall say
that  $G$ {\it has the property $C_\pi$} ($G\in C_\pi$). If we have further
that each  $\pi$-subgroup of $G$ is contained in a  $\pi$-Hall
subgroup, then we shall say that $G$ {\it has the property $D_\pi$} ($G\in
D_\pi$). A group with the property $E_\pi$ ($C_\pi$, $D_\pi$) we shall call
also an $E_\pi$-
(respectively a $C_\pi$-, a $D_\pi$-) {\it group}. The expression
(mod CFSG) means that the corresponding result is proved by using the
classification of finite simple groups.

Assume that a set $\pi$ is fixed. It is proved that the class of all
$D_\pi$-groups is closed under homomorphic images, normal subgroups (mod CFSG,
\cite[Theorem~7.7]{RevVdoDpiFinal}) and extensions (mod
CFSG, \cite[Theorem~7.7]{RevVdoDpiFinal}). The class of $E_\pi$-groups is also
known to be closed under normal subgroups and homomorphic images (see Lemma
\ref{base}(1)), but, in general, it is not closed under extensions (see
Example 1). The class of $C_\pi$-groups is closed under homomorphic images (mod
CFSG, see Lemma \ref{quot}) and extensions (see Lemma \ref{cpiext}), but, in
general, is not closed under normal subgroups (see Example 2). Nevertheless,
while proving statements about Hall properties one need to know in which cases
an extension of an $E_\pi$-group by an $E_\pi$-group has the property $E_\pi$,
and a normal subgroup of a $C_\pi$-group has the property $C_\pi$. In the
present paper we prove by using the classification of finite simple groups

\begin{Theo}\label{MainTheorem} {\em (mod CFSG)}
Let $\pi$ be a set of primes, $H$ be a $\pi$-Hall, and $A$ be a normal subgroups
of a $C_\pi$-group $G$. Then~${HA\in C_\pi}$.
\end{Theo}

Since every $\pi'$-group has the property $C_\pi$ and the class of
$C_\pi$-groups is closed under extensions, the following statement is immediate
from Theorem \ref{MainTheorem} (cf. Example 2).

\begin{Cor}\label{FirstCorollary}
Let $\pi$ be a set of primes and $A$ be a normal subgroup of
$G$, which index in $G$ is a $\pi'$-number. Then $G\in C_\pi$ if and only
if~${A\in C_\pi}$.
\end{Cor}

The proof of Theorem \ref{MainTheorem} is based on the following result about
the number of classes of $\pi$-Hall subgroups in finite simple groups.

\begin{Theo} \label{Simple} {\em (mod CFSG)}
Let $\pi$ be a set of primes and  $S$ a finite simple group. Then there exist
at most $4$ classes of conjugate $\pi$-Hall subgroups of $S$. More precisely,
the following statements hold.

$(1)$ If $2\not\in\pi$, then $S$ has at most one class of conjugate
$\pi$-Hall subgroups.

$(2)$ If $2\in\pi$, and $3\not\in\pi$, then $S$ has at most two classes of
conjugate $\pi$-Hall subgroups.

$(3)$ If $2,3\in\pi$, then $S$ has at most four classes of
conjugate $\pi$-Hall subgroups.
\end{Theo}

\begin{Cor}\label{SimpleCor}
Let $\pi$ be a set of primes and $S$ a finite simple $E_\pi$-group. Then if a
positive integer $k$ is not greater than the number of classes of conjugate
$\pi$-Hall subgroups of $S$, then $k$ is a $\pi$-num\-ber.
\end{Cor}

In view of Theorem \ref{MainTheorem} note that the authors do not know
any counter example to the following conjecture.

\begin{Conj}
Let $\pi$ be a set of primes,  $G$ be a finite
$C_\pi$-group, and $A$ be its (not necessary normal) subgroup, containing a
$\pi$-Hall subgroup of $G$. Then $A\in C_\pi$.
\end{Conj}

In the conjecture the condition that $A$ contains a $\pi$-Hall
subgroup of  $G$ cannot be weeken by substituting it with condition
that the index of $A$ is a $\pi'$-number. Indeed, consider
$B_3(q)\simeq \mathrm{P}\Omega_7(q)$, where $q-1$ is divisible by
$12$ and is not divisible by  $24$. In view of
\cite[Lemma~6.2]{RevVdoDpiFinal} the group $\mathrm{P}\Omega_7(q)$
is a $C_{\{2,3\}}$-group and its $\{2,3\}$-Hall subgroup is
contained in a monomial subgroup. On the other hand the group
$\Omega_7(2)$ is known to be embeddable into $\mathrm{P}\Omega_7(q)$
and, under above mentioned conditions on  $q$, it index is not
divisible by $2$ and by $3$. Although $\Omega_7(2)$ does not contain
$\{2,3\}$-Hall subgroups, i.~e., it is not even a
$E_{\{2,3\}}$-group.

\section{Notations and preliminary results}

By $\pi$ we always denote a set of primes, and the term group always means
a finite group.

For a group  $G$, a $G$-class of $\pi$-Hall subgroups is a class of conjugate
$\pi$-Hall subgroups of $G$. Let $A$ be a subnormal subgroup of a
$E_\pi$-group $G$. A subgroup of $A$ of type $H\cap A$, where $H$ is a
$\pi$-Hall subgroup of $G$, we shall call a {\it
$G$-in\-du\-ced $\pi$-Hall subgroup} of $A$. Thus a set $\{(H\cap A)^x\mid x\in
A\}$, where $H$ is a $\pi$-Hall subgroup of  $G$ is called an {\it
$A$-class of $G$-in\-du\-ced $\pi$-Hall subgroups}. By $k_\pi^G(A)$ we denote
the number of all $A$-clas\-ses of $G$-in\-du\-ced $\pi$-Hall subgroups. Assume
also that  $k_\pi(G)=k_\pi^G(G)$ is the number of conjugacy classes of
$\pi$-Hall subgroups of $G$. It is clear that~${k_\pi^G(A)\le k_\pi (A)}$.

Let $A=A_1\times\dots\times A_s$ and for each $i=1,\dots,s$ by ${\mathscr
K}_i$ we denote an $A_i$-class of $\pi$-Hall subgroups of $A_i$. The set
$${\mathscr K}_1\times \dots\times {\mathscr K}_s=\{\langle H_1,\dots,
H_s\rangle\simeq H_1\times\ldots\times H_s\mid H_i\in {\mathscr
K}_i, i=1,\dots,s\}$$ is called a {\it product of classes}
${\mathscr K}_1, \dots, {\mathscr K}_s$. Clearly ${\mathscr
K}_1\times \dots\times {\mathscr K}_s$ is an $A$-class of $\pi$-Hall
subgroups of $A$. It is also clear that if $A$ is a normal subgroup
of  $G$, then each $A$-classes of $G$-in\-du\-ced $\pi$-Hall
subgroups is a product of some  $A_1$-, $\dots$, $A_s$- classes of
$G$-in\-du\-ced $\pi$-Hall subgroups. The reverse statement, in
general, is not true.

The following statements are known and their proof does not use the
classification of finite simple groups.

\begin{Lemma} \label{base}
Let $A$ be a normal subgroup of $G$. Then the following statements hold.

{\rm (1)} If $H$ is a $\pi$-Hall subgroup of $G$, then
$H \cap A$ is a $\pi$-Hall subgroup of $A$, and $HA/A$ is a
$\pi$-Hall subgroup of~$G/A$.

{\rm (2)} If all factors of a subnormal series of $G$ are either $\pi$- or
$\pi'$- groups, then $G\in D_{\pi}$.
\end{Lemma}

\begin{Lemma}\label{cpiext} {\em (S.A.Chunihin, see also \cite[Theorems~C1
and~C2]{Hall})}
Let $A$ be a normal subgroup of $G$. If $A$ and $G/A$ has the property
$C_{\pi}$, then~${G\in C_{\pi}}$. 
\end{Lemma}

\noindent{\bf Example 1.} {Let $\pi=\{2,3\}$. Let $G=\GL_3(2)=\SL_3(2)$ be a
group of order $168=2^3\cdot 3\cdot 7$. From
\cite[Theorem  1.2]{RevHallp} or \cite{Atlas} it follows that $G$ has exactly
two classes of $\pi$-Hall subgroups with representatives
$$
\left(
\begin{array}{c@{}c}
\fbox{
$\begin{array}{c}
\\
\!\!\!
\GL_2(2)
\!\!\! \\
\\
\end{array}$
}
& *\\
0 &\fbox{1}
\end{array}
\right)\text{ and }
\left(
\begin{array}{c@{}c}
\fbox{1}& *\\
 0&\fbox{$\begin{array}{c}
\\
\!
\GL_2(2)
\!
\\
\\
\end{array}$}
\end{array}
\right).
$$
The first one consists of  line stabilizers in the natural
representation of $G$, and the second one consists of plain
stabilizers. The map $\iota : x\in G\mapsto (x^t)^{-1}$, where $x^t$
means the transposed matrix to  $x$, is an automorphism of order $2$
of  $G$. It interchange classes of $\pi$-Hall subgroups. If
$\widehat{G}=G\leftthreetimes\langle\iota\rangle$ is a natural
semidirect product, then $N_{\widehat{G}}(H)=N_G(H)=H$ for each
$\pi$-Hall subgroup $H$ of $G$, since each element from
$\widehat{G}\setminus G$ interchanges classes of $\pi$-Hall
subgroups, like  $\iota$ does, and $H$ is maximal. The group
$\widehat{G}$ has order $2^4\cdot3\cdot7$ and has no the property
$E_\pi$. Indeed, if there would exist a $\pi$-Hall subgroup
$\widehat{H}$, then by Lemma \ref{base}(1) there would also exist a
$\pi$-Hall subgroup $H$ of $G$, with $H=\widehat{H}\cap G$. But this
implies  $\widehat{H}\le N_{\widehat{G}}(H)=H$, and so
$|\widehat{H}|\le 2^3\cdot 3$, that contradicts the condition that
$\widehat{H}$ is a $\pi$-Hall subgroup.

\noindent{\bf Example 2.} {Let $\pi=\{2,3\}$. Let $G=\GL_5(2)=\SL_5(2)$ be a
group of order $99999360=2^{10}\cdot 3^2\cdot 5\cdot 7 \cdot
31$. Let $\iota : x\in G\mapsto (x^t)^{-1}$ and $\widehat{G}=G\leftthreetimes
\langle\iota\rangle$ be a natural semidirect product. From
\cite[Theorem~1.2]{RevHallp} it follows that there exist $\pi$-Hall subgroups
of $G$, and every such a subgroup is a stabilizer of a series of subspaces
$V=V_0<V_1<V_2<V_3=V$, where $V$ is a natural module of $G$, and $\dim
V_k/V_{k-1}\in \{1,2\}$ for every $k=1,2,3$. Therefore, there are
three conjugacy classes of $\pi$-Hall subgroups of $G$ with representatives
$$
H_1=\left(
\begin{array}{c@{}c@{}c}
\fbox{$\begin{array}{c}
\\
\!
\GL_2(2)
\!
\\
\\
\end{array}$}& &*\\
 &\fbox{1}& \\
0& &
\fbox{$\begin{array}{c}
\\
\!
\GL_2(2)
\!
\\
\\
\end{array}$}
\end{array}
\right),
$$
$$H_2=
{\left(
\begin{array}{c@{}c@{}c}
\fbox{1}& &*\\
 &\fbox{$\begin{array}{c}
\\
\!
\GL_2(2)
\!
\\
\\
\end{array}$}& \\
0& &
\fbox{$\begin{array}{c}
\\
\!
\GL_2(2)
\!
\\
\\
\end{array}$}
\end{array}
\right),
\text{ and }
H_3=\left(
\begin{array}{c@{}c@{}c}
\fbox{$\begin{array}{c}
\\
\!
\GL_2(2)
\!
\\
\\
\end{array}$}& &*\\
 &\fbox{$\begin{array}{c}
\\
\!
\GL_2(2)
\!
\\
\\
\end{array}$}& \\
0& &
\fbox{1}
\end{array}
\right).}
$$
Note that  $N_G(H_k)=H_k$, $k=1,2,3$, since  $H_k$ is parabolic. The class
containing $H_1$ is $\iota$-invariant. So Frattini argument implies that
$\widehat{G}=GN_{\widehat{G}}(H_1)$, whence $|N_{\widehat{G}}(H_1):N_G(H_1)|=2$
and so $N_{\widehat{G}}(H_1)$ is a $\pi$-Hall subgroup of  $\widehat{G}$.
$\iota$ interchanges classes containing  $H_2$ and $H_3$. So, as in the
previous example, these subgroups are not contained in $\pi$-Hall subgroups of
$\widehat{G}$. Thus  $\widehat{G}$ contains exactly one class of $\pi$-Hall
subgroups, therefore has the property $C_\pi$, while its normal subgroup  $G$
has no this property.}

\begin{Lemma} \label{Vedernik}
{\em (\cite[Lemma 3]{Ved}, mod CFSG)} Let $A$ be a normal subgroup of an
$E_\pi$-group $G$, $B$ be a subgroup of  $G$ such that
$A\leq B$ and $B/A\in E_\pi$. Then  $B\in E_\pi$.
\end{Lemma}

It is immediate from Lemma \ref{Vedernik} that the property $C_\pi$
is preserved under homomorphisms and we shall give the proof here
for completeness. Note that since Lemma \ref{Vedernik} is proved by
using the classification of finite simple groups, then Lemma
\ref{quot} is also proved by using the classification of finite
simple groups.

\begin{Lemma} \label{quot} {\em (mod CFSG)}
Let $A$ be a normal subgroup of a $C_\pi$-group $G$. Then $G/A\in
C_\pi$.
\end{Lemma}

\begin{proof}
Let $G/A=\overline{G}$. Since all $\pi$-Hall subgroups of $G$ are
conjugate, it is enough to prove that for every $\pi$-Hall subgroup
$\overline{K}$ of $\overline{G}$ there exists a $\pi$-Hall subgroup
$U$ of $G$ such that $UA/A=\overline{K}$. The existence of such a
subgroup $U$ follows from Lemma \ref{Vedernik}. Indeed, for $K$
being the complete preimage of $\overline{K}$, the group
$K/A=\overline{K}$ has the property $E_\pi$. So $K$ has the property
$E_\pi$, and its $\pi$-Hall subgroup  $U$ is a $\pi$-Hall subgroup
of $G$. By construction,~${UA/A=\overline{K}}$.
\end{proof}

\begin{Lemma} \label{Norm}
Let $A$ be a normal and  $H$ be a $\pi$-Hall subgroups of a
$C_\pi$-group $G$. Then each of groups $N_G(HA)$ and $N_G(H\cap A)$
has the property~$C_\pi$.
\end{Lemma}

\begin{proof}
Let $K$ be a $\pi$-Hall subgroup of $N_G(HA)$. Since $HA\unlhd
N_G(HA)$ and $|N_G(HA):HA|$ is a $\pi'$-number, then $K\leq HA$ and
$KA=HA$. If $x\in G$ is such that $K=H^x$, then $(HA)^x=H^xA=KA=HA$
and thus $x\in N_G(HA)$. Therefore $N_G(HA)\in C_\pi$.

Assume now that $K$ is a  $\pi$-Hall subgroup of $N_G(H\cap A)$.
Then $K(H\cap A)=K$ and $K\cap A=H\cap A$. If $x\in G$ is such that
$K=H^x$, then $(H\cap A)^x=H^x\cap A=K\cap A=H\cap A$ and thus $x\in
N_G(H\cap A)$. Therefore $N_G(H\cap A)\in C_\pi$.
\end{proof}

\begin{Lemma} \label{ind=inv}
Let $A$ be a normal, $H$ be a $\pi$-Hall subgroup of a group $G$ and
$HAC_G(A)\unlhd G$ (this condition is satisfied, if $HA\unlhd G$).
Then an $A$-class of $\pi$-Hall subgroups is a class of
$G$-in\-du\-ced  $\pi$-Hall subgroups if and only if it is
$H$-in\-va\-ri\-ant.
\end{Lemma}

\begin{proof}
If $K$ is a $\pi$-Hall subgroup of $G$, then $K\leq HAC_G(A)$ holds
and so $KAC_G(A)=HAC_G(A)$. Since the $A$-class $\{(K\cap A)^x\mid
x\in A\}$ is $K$-in\-va\-ri\-ant, it is invariant under
$HAC_G(A)=KAC_G(A)$ and hence under~$H$.

Now we prove the inverse statement. Without lost of generality we
may assume that  $G=HA$. Let $U$ be a $\pi$-Hall subgroup of $A$
such that $\{U^x\mid x\in A\}$ is $H$-in\-va\-ri\-ant. Then Frattini
argument implies $H\leq N_G(U)A$ and so $G=N_G(U)A$. Note that
$$N_G(U)/N_A(U)=N_G(U)/N_G(U)\cap A\simeq N_G(U)A/A=G/A=HA/A $$ is a
$\pi$-group. Thus $N_G(U)$ has a normal series
$$N_G(U)\geq N_A(U)\geq U\geq 1,$$
each  factor of this series is either a $\pi$- or a
$\pi'$- group. By Lemma \ref{base}(2) $N_G(U)$ has the property $D_\pi$. Let $K$
be its
$\pi$-Hall subgroup. Clearly $U\leq K$. More over $K$ is a
$\pi$-Hall subgroup of $G$, since
$$|K|=|N_G(U)|_\pi=|N_G(U)/N_A(U)|_\pi|N_A(U)|_\pi=|HA/A|_\pi|U|=|HA/A||H\cap
A|=|H|.$$ Therefore, $U=K\cap A$ is a $G$-in\-du\-ced $\pi$-Hall
subgroup.
\end{proof}

\begin{Lemma} \label{EqualityG=HA}
Let $A$ be a normal,  $H$ be a $\pi$-Hall subgroups of a
$C_\pi$-group $G$ and $HA\unlhd G$. Then $k_\pi^G(A)=k_\pi^{HA}(A)$.
\end{Lemma}

\begin{proof}
Since $HA$ is a normal subgroup of  $G$, each $\pi$-Hall subgroup of
$G$ is contained in $HA$. By Lemma \ref{Norm} $HA\in C_\pi$, and so
the equality ${k_\pi^G(A)=k_\pi^{HA}(A)}$ holds.
\end{proof}

\begin{Lemma}\label{crit}
Let $A$ be a normal, $H$ be a $\pi$-Hall subgroup of a $E_\pi$-group
$G$ and $HA\unlhd G$. Then the following statements are equivalent.

$(1)$ $k_\pi^G(A)=1$.

$(2)$ $HA\in C_\pi$.

$(3)$ Every two $\pi$-Hall subgroups of $G$ are conjugate by an
element of~$A$.
\end{Lemma}

\begin{proof}
$(1)\Rightarrow (2)$. If $K$ is a $\pi$-Hall subgroup of  $HA$, then
by (1) groups $H\cap A$ and $K\cap A$ are conjugate in $A$. We may
assume that $H\cap A=K\cap A$. Then  $H$ and $K$ are contained in
$N_{HA}(H\cap A)$. By Frattini argument we have $HA=N_{HA}(H\cap
A)A$. So
$$N_{HA}(H\cap A)/N_A(H\cap A)=N_{HA}(H\cap A)/N_{HA}(H\cap A)\cap A\simeq
N_{HA}(H\cap A)A/A=HA/A$$ is a $\pi$-group. Thus  $N_{HA}(H\cap A)$
has a normal series
$$N_{HA}(H\cap A)\geq N_A(H\cap A)\geq H\cap A\geq 1,$$
each factor of this series is either a $\pi$- or a $\pi'$- group,
and, by Lemma \ref{base}(2) has the property $D_\pi$. In particular,  $H$ and
$K$ are
conjugate in $N_{HA}(H\cap A)$.

 $(2)\Rightarrow (3)$ and  $(3)\Rightarrow (1)$ are evident.
\end{proof}

 \begin{Lemma} \label{TrivAct}
Let $A=A_1\times\dots\times A_s$ be a normal and $H$ be a $\pi$-Hall
subgroups of $G$. Assume also that subgroups $A_1,\ldots,A_s$ are
normal in $G$ and $G=HA C_G(A)$. Then $k_\pi^G(A)=k_\pi^G(A_1)\dots
k_\pi^G(A_s)$.
\end{Lemma}

\begin{proof}
Two $\pi$-Hall subgroups $P$ and $Q$ of $A$ are conjugate in $A$ if
and only if $\pi$-Hall subgroups $P\cap A_i$ and $Q\cap A_i$ of
$A_i$ are conjugate in $A_i$ for each $i=1,\dots,s$. To prove the
lemma it is enough to show that the product of every $A_1$-,
$\dots$, $A_s$- classes of $G$-in\-du\-ced $\pi$-Hall subgroups is
an $A$-class of again $G$-in\-du\-ced $\pi$-Hall subgroups. Let
$U_1,\dots,U_s$ be $G$-in\-du\-ced $\pi$-Hall subgroups of
$A_1,\dots,A_s$ respectively. We shall show that
$$U=\langle U_1,\dots,U_s\rangle=U_1\times\ldots\times U_s$$ is a $G$-in\-du\-ced
$\pi$-Hall subgroup of $A$. By Lemma \ref{ind=inv} it is enough to
show that for each  $h\in H$ there exists $a\in A$ with $U^h=U^a$.
Since $U_i=K_i\cap A_i$ for a suitable $\pi$-Hall subgroup $K_i$ of
$G$, then $\{U_i^{x_i}\mid x_i\in A_i\}$ is invariant under $K_i$
and so under $K_iA=HA$. In particular, $U_i^h=U_i^{a_i}$ for some
$a_i\in A_i$. Thus
$$U^h=U_1^h\times\ldots\times U_s^h=
U_1^{a_1}\times\ldots\times U_s^{a_s}=U_1^a\times\ldots\times
U_s^a=U^a,$$ where $a=a_1\ldots a_s\in A$.
\end{proof}

\begin{Lemma} \label{TransAct}
Let $A=A_1\times\dots\times A_s$ be a normal, $H$ be a $\pi$-Hall
subgroups of $G$. Assume that $G$ acts transitively by conjugation
on the set $\{A_1,\dots,A_s\}$ and $G=HA C_G(A)$. Then
$k_\pi^G(A)=k_\pi^G(A_1)=\dots= k_\pi^G(A_s)$.

\end{Lemma}

\begin{proof}
We shall show that if $x,y\in G$ are in the same coset of  $G$ by
$N_G(A_1)$, then for each $G$-in\-du\-ced  $\pi$-Hall subgroup $U_1$
of $A_1$ subgroups $U_1^x$ and $U_1^y$ are conjugate in
$A_i=A_1^x=A_1^y$. It is enough to show that subgroups  $U_1$ and
$U_1^t$, where $t=xy^{-1}\in N_G(A_1)$, are conjugate in $A_1$. Let
$t=ach$, $a\in A$, $c\in C_G(A)$, $h\in H$. Since subgroups $U_1$
and $U_1^{ac}$ are conjugate in $A_1$, we need to show that $U_1$
and $U_1^h$ are conjugate in $A_1$. Since $ac\in N_G(A_1)$, then $h$
normalizes also the subgroup $A_1$. Let $U_1=U\cap A_1$ for some
$G$-in\-du\-ced $\pi$-Hall subgroup $U$ of $A$. Let ${\mathscr K}$
be an $A$-class, containing $U$, and let ${\mathscr K}={\mathscr
K}_1\times \dots\times {\mathscr K}_s$, where ${\mathscr K}_i$ is an
$A_i$-class of $G$-in\-du\-ced $\pi$-Hall subgroups. Clearly $U_1\in
{\mathscr K}_1$. Since, according to Lemma \ref{ind=inv}, the class
${\mathscr K}$ is $H$-invariant, then $H$ acts on the set
$\{{\mathscr K}_1, \dots, {\mathscr K}_s\}$. An element $h$
normalizes $A_1$, so stabilizes the $A_1$-class  ${\mathscr K}_1$.
In particular, subgroups $U_1$ and $U_1^h$ are in ${\mathscr K}_1$
and so are conjugate in~$A_1$.

Let $f\in H$ and $A_1^f=A_i$. For an $A_1$-class of $\pi$-Hall subgroups
${\mathscr K}_1$ define an $A_i$-class ${\mathscr K}_1^{f}$, assuming
${\mathscr
K}_1^{f}=\{U_1^{f}\mid U_1\in {\mathscr K}_1\}$. Clearly ${\mathscr
K}_1^{f}$ is an $A_i$-class of $G$-in\-du\-ced $\pi$-Hall subgroups.

Let $h_1,\dots, h_s$ be a right transversal of $H$ by $N_H(A_1)$, with
$h_1=1$. Up to the renumberring we may assume that
$A_i=A_1^{h_i}$. Since $AC_G(A)\leq N_G(A_i)$ for every
$i=1,\dots,s$, then elements $h_1,\dots, h_s$ forms a right transversal of $G$
by $N_G(A_1)$. We shall show that the map
$${\mathscr K}_1\mapsto {\mathscr K}_1^{h_1}\times\dots\times{\mathscr
K}_1^{h_s}$$ is a bijection between the set of $A_1$-clas\-ses of
$G$-in\-du\-ced $\pi$-Hall subgroups and the set of $A$-clas\-ses of
$G$-in\-du\-ced $\pi$-Hall subgroups. Since ${\mathscr
K}_1^{h_1}\times\dots\times{\mathscr K}_1^{h_s}$ is
$H$-in\-va\-ri\-ant, then by Lemma \ref{ind=inv} its elements are 
$G$-in\-du\-ced $\pi$-Hall subgroups. It is also clear that the map is
injective. Let ${\mathscr K}={\mathscr K}_1\times
\dots\times {\mathscr K}_s$ be an $A$-class of $G$-in\-du-ced $\pi$-Hall
subgroups. To prove the surjectivity it is enough to show that  ${\mathscr
K}_i={\mathscr K}_1^{h_i}$ for every $i=1,\dots, s$. Since $G$ acts
transitively on $\{A_1,\dots,A_s\}$, then there exists an element $g\in G$ such
that $A_1^g=A_i$. Let $g=act$, where $a\in A$, $c\in C_G(A)$, and $t\in H$.
Then $A_1^t=A_i$ and $t\in N_H(A_1)h_i$. As we have proved,  ${\mathscr
K}_1^t={\mathscr K}_1^{h_i}$. By Lemma \ref{ind=inv} the class 
${\mathscr K}$ is $H$-in\-va\-ri\-ant, therefore ${\mathscr
K}_i={\mathscr
K}_1^t={\mathscr K}_1^{h_i}$ and ${\mathscr K}={\mathscr
K}_1^{h_1}\times\dots\times{\mathscr K}_1^{h_s}.$
\end{proof}

\section{Proof of Theorem \ref{Simple}}

Statement (1) of the theorem follows from \cite[Theorem~А]{GroConjOddOrder}.
Statement (2) of the theorem follows from \cite[Lemma~5.1 and
Corollary~5.4]{RevVdoDpiFinal} Thus we need to prove  statement (3) of the
theorem.

Let $\pi$ be a set of primes such that $2,3\in\pi$. Consider all known
simple groups separately.

If $S\simeq \mathrm{Alt}_n$, then the statement of the lemma follows from
\cite[Theorem~4.3]{RevVdoDpiFinal}.

In sporadic groups all proper $\pi$-Hall subgroups with $2,3\in\pi$ are given in
\cite[Table~2]{RevDpiOneClass} (note that there is an evident missprint
in the table, $\{2,3,7\}$-Hall subgroup of $J_1$ has the structure
$2^3:7:3$). By using \cite{Atlas}, it is easy to see that in all
cases~${k_\pi(S)\le2}$.

Let $S$ be a finite group of Lie type over a field of characteristic  $p\in\pi$
and $H$ be a $\pi$-Hall subgroup of $S$. By \cite[Theorem~1.2]{RevHallp} one of
the following cases 1--4 holds.

Case 1. $H=S$. Clearly $k_\pi(S)=1$.

Case 2. $H$ is contained in a Borel subgroup of $S$. In this case, since Borel 
subgroups are conjugate and solvable, we have that~${k_\pi(S)=1}$.

Case 3.
$S\simeq D_n^\varepsilon(2^t)\simeq \mathrm{P}\Omega_{2n}^\varepsilon(q)$,
$\varepsilon\in\{+,-\}$, $H$ is a parabolic subgroup with a Levi factor
isomorphic to $D_{n-1}^\varepsilon(2^t)$. Since all such parabolic subgroups
are conjugate,~${k_\pi(S)=1}$.

Case 4.
$S\simeq {\rm PSL}(V)$ and $H$ is an image in ${\rm PSL}(V)$ of a parabolic
subgroup of $\SL(V)$, stabilizing a series of subspaces
$0=V_0<V_1<\ldots<V_s=V$. If $n_i$ is equal to  $\dim
V_i/V_{i-1}$, and $n$ is equal to
$\dim V$, then one of the following subspaces holds

(4.1) $n$ is an odd prime, $s=2$, $\{n_1,n_2\}=\{1,n-1\}$,
$k_\pi(S)=2$;

(4.3)
$n=4$, $s=2$, $n_1=n_2=2$, $k_\pi(S)=1$;

(4.4) $n=5$, $s=2$ $\{n_1,n_2\}=\{2,3\}$, $k_\pi(S)=2$;

(4.5) $n=5$, $s=3$, $\{n_1,n_2,n_3\}=\{1,2,2\}$, $k_\pi(S)=3$;

 (4.6) $n=7$, $s=2$, $\{n_1,n_2\}=\{3,4\}$, $k_\pi(S)=2$;

(4.7) $n=8$, $s=2$, $n_1=n_2=4$, $k_\pi(S)=1$;

(4.8) $n=11$, $s=2$, $\{n_1,n_2\}=\{5,6\}$, $k_\pi(S)=2$.

Assume at last that $S$ is a finite group of Lie type over a field of
characteristic $p$, $p\not\in\pi$ and $H$ is a $\pi$-Hall subgroup of~$S$.

If $S\simeq A_n^\varepsilon(q)$, where $n\ge2$, $\varepsilon\in\{+,-\}$,
$A_n^+(q)=A_n(q)$ and
$A_n^-(q)={}^2A_n(q)$, then by \cite[Lemma~6.1]{RevVdoDpiFinal} we have
$k_\pi(S)=1$. Assume that $S\simeq A_1(q)$. If 
$\pi\cap\pi(S)\not= \{2,3\}$ and $\pi\cap\pi(S)\not= \{2,3,5\}$,
then by \cite[Lemma~6.1]{RevVdoDpiFinal} we have $k_\pi(S)=1$.

Assume that $\pi\cap\pi(S)=\{2,3\}$. Then by \cite[Lemma~6.1]{RevVdoDpiFinal},
either $H$ is contained in the normalizer of a maximal torus (and all 
such subgroups are conjugate), or $H$ is isomorphic to
$\mathrm{Alt}_4$ or $\mathrm{Sym}_4$ and $H$ is a homomorphic image of a
primitive absolutely irreducible solvable subgroup of
$\SL_2(q)$. From \cite[\S~21, Theorem~6]{SupMat}it follows that all such
subgroups $H$ are conjugate in ${\rm PGL}_2(q)$, hence in
${\rm PSL}_2(q)$ there exist at most $2$ classes of such subgroups. Thus 
in this case~${k_\pi(S)\le 3}$.

Assume now that $\pi\cap\pi(S)=\{2,3,5\}$. By
\cite[Lemma~6.1]{RevVdoDpiFinal} it follows that in this case either $H$ is
contained in the normalizer of a maximal torus (and all such subgroups are
conjugate), or $H$ is isomorphic to ${\rm PSL}_2(5)=\mathrm{Alt}_5$. More over
since $p\not\in\pi$, then $p$ does not divide $\vert \SL_2(5)\vert$, in
particular, $p\not=2,3,5$. By \cite{Atlas} or \cite{GAP} it follows that
$\SL_2(5)$, being the preimage of $H$ in $\SL_2(q)$, there exist $2$ exact
irreducible complex representation of degree $2$, conjugated by an outer
automorphism of $\SL_2(5)$, while  $\mathrm{Alt}_5$ has no representations of
degree $2$. By \cite[Theorem~15.3 and Corollary~9.7]{Isaacs} and the fact that 
$\SL_2(5)\leq \GL_2(q)$, there exist precisely two nonequivalent representations
of degree $2$ over the field  $\mathrm{F}_q$ of $\SL_2(5)$, conjugated by an
outer automorphism. Therefore, all subgroups  $H$, isomorphic to ${\rm
PSL}_2(5)$, are conjugate in ${\rm PGL}_2(q)$, ans so in ${\rm
PSL}_2(q)$ there exists at most $2$ classes of such subgroups. Thus, in this
case $k_\pi(S)\le 3$.

Let $S$ be a group of Lie type
distinct from $A_n^\varepsilon(q)$, $H$ is its $\pi$-Hall subgroup, and
$2,3\in\pi$, $p\not\in\pi$ (in particular, $p$ does not divide $\vert H\vert$
and $\vert 2.H\vert$). If $S$ is not isomorphic to one of the groups
$\mathrm{P}\Omega_7(q)$, $\mathrm{P}\Omega_8^+(q)$ or $\mathrm{P}\Omega_9(q)$,
then by \cite[Lemmas~6.2--6.10]{RevVdoDpiFinal}, $H$ is contained in the
normalizer of a maximal torus and all such $\pi$-Hall subgroups are conjugate.
If $S$ is isomorphic to one of the groups
$\mathrm{P}\Omega_7(q)$,
$\mathrm{P}\Omega_8^+(q)$ or $\mathrm{P}\Omega_9(q)$, then by
\cite[Lemmas~6.2 and~6.4]{RevVdoDpiFinal} either $H$ is solvable and
contained in the normalizer of a maximal torus (in this case all such 
$\pi$-Hall subgroups are conjugate and the order of a Sylow $7$-subgroup
is greater than $7^2$), or
$\pi\cap \pi(S)=\{2,3,5,7\}$ and the pair $(S,H)$ lies in the following list:
$(\mathrm{P}\Omega_7(q),\Omega_7(2))$,
$(\mathrm{P}\Omega_8^+(q),\Omega_8^+(2))$,
$(\mathrm{P}\Omega_9(q),{2.\Omega_8(2)^+\!:\!2})$. Note that if $H$ is
nonsolvable, then the order of its Sylow  $7$-subgroup is not greater than
$7^2$, so there does not exists a $\pi$-Hall subgroup contained in the
normalizer of a maximal torus. By \cite[Lemma~1.7.1]{KleidO8} the number of
conjugacy classes of subgroups of
${\rm GO}_7(q)$, isomorphic to $\Omega_7(2)$, and subgroups of ${\rm
GO}_8^+(q)$, isomorphic to $2.\Omega_8^+(2)$, is not greater than the number of
irreducible representation of these subgroups of degrees $7$ and $8$
respectively. By using ordinary characters tables of $\Omega_7(2)$ and
$2.\Omega_8^+(2)$, given in \cite{Atlas} or \cite{GAP}, we obtain that both
these groups has precisely one irreducible complex representation of degree 
$7$ or $8$ respectively, while $2.\Omega_7(2)$ and $\Omega_8^+(2)$ have no
complex representations of degrees $7$ and $8$ respectively. By
\cite[Theorem~15.3 and Corollary~9.7]{Isaacs}, this statement also holds for 
representations over $\mathrm{F}_q$. It follows the inequalities
$k_\pi(S)\le2$ if $S\simeq\mathrm{P}\Omega_7(q)$ and $k_\pi(S)\le4$ if
$S\simeq\mathrm{P}\Omega_8^+(q)$. In the case 
$(S,H)=(\mathrm{P}\Omega_9(q),2.\Omega_8(2)^+\!:\!2)$ the group 
$H$ is contained in the centralizer of an involution an by
\cite[Proposition~11]{Ca2} all such centralizers are conjugate, so in this case
$k_\pi(S)\le 4$ also.

\section{Proof of Theorem \ref{MainTheorem}}

Assume that the theorem is not true and $G$ is a counter example of minimal
order. Then $G$ contains a $\pi$-Hall subgroup $H$ and a normal subgroup $A$
such that $HA$ has no the property $C_\pi$. Choose from such subgroups $A$ the
minimal by inclusion. Let $K$ be a $\pi$-Hall subgroup of $HA$, that is not
conjugate to $H$ in $HA$. The process of eliminating of $G$ we divide into
several steps.
\smallskip

Clearly

\smallskip
{\it $(1)$ $HA=KA$.}
\smallskip

{\it $(2)$ $A$ is a minimal normal subgroup of $G$.}
\smallskip

Otherwise let $M$ be a nontrivial normal subgroup of $G$, that is properly
contained in $A$. Let $\overline{G}=G/M$, and for each subgroup $B$ of $G$ by
$\overline{B}$ we denote the group $BM/M$. By Lemma \ref{quot} the group 
$\overline{G}$ has the property $C_\pi$, $\overline{H}$ and $\overline{K}$ are
its $\pi$-Hall subgroups, $\overline{A}$ is a normal subgroup, 
$\overline{H}\overline{A}=\overline{K}\overline{A}$ and $|\overline{G}|<|G|$.
By the minimality of the counter example, $G$ the group 
$\overline{H}\overline{A}$ has the property $C_\pi$. So subgroups
$\overline{H}$ and  $\overline{K}$ are conjugate by an element of 
$\overline{A}$. This means that subgroups $HM$ and $KM$ are conjugate by an
element of $A$. Without lost of generality we may assume that  $HM=KM$. In view
of the choice of $A$, the group  $HM$ has the property $C_\pi$. But this means
that  $H$ and $K$ are conjugate by an element of $M\leq A$, a contradiction.

\smallskip
{\it $(3)$ $A\not\in C_\pi$. In particular, $A$ is not solvable.}
\smallskip

Otherwise by Lemma \ref{cpiext} $HA$ would have the property $C_\pi$ as an
extension of a 
$C_\pi$-group by a  $\pi$-group.

\smallskip
{\it $(4)$ $HA$ is a normal subgroup of $G$.}
\smallskip

Otherwise $N_G(HA)$ is a proper subgroup of $G$ and by Lemma \ref{Norm}
$N_G(HA)\in C_\pi$. Since $G$ is a counter example of minimal order,   $HA\in
C_\pi$, a contradiction.
\smallskip

By $(2)$ and $(3)$

\smallskip
{\it $(5)$  $A$ is a direct product of simple non-Abelian groups
$S_1,\dots,S_m$. The group $G$ acts transitively by conjugation on
$\Omega=\{S_1,\dots,S_m\}$.}
\smallskip

Let $\Delta_1,\dots\Delta_s$ be orbits of $HA$ on the set
$\Omega$, and let $T_j=\langle\Delta_j\rangle$ for each 
$j=1,\dots,s$. By $(4)$ and $(5)$

\smallskip
{\it $(6)$ $G$ acts transitively by conjugation on 
$\{T_1,\dots T_s\}$. The subgroup $A$ is a direct product of
$T_1,\dots T_s$, and each of $T_i$ is normal in~$HA$. }
\smallskip

Let $S\in\Omega$ and $T$ is the subgroup, generated by the orbit from
$\Delta_1,\dots\Delta_s$, that contains $S$. By Lemmas лемм \ref{crit} and 
\ref{TransAct} it follows that

\smallskip
{\it $(7)$ $k_\pi^G(T)=k_\pi^G(S)$.}
\smallskip

By $(7)$ and Lemmas \ref{crit} and \ref{TrivAct} it follows

\smallskip
{\it $(8)$ $k_\pi^G(A)=(k_\pi^G(T))^s=(k_\pi^G(S))^s$.}
\smallskip

Since $1\le k_\pi^G(S)\le k_\pi(S)$, from $(8)$ and by Corollary
\ref{SimpleCor} it follows that

\smallskip
{\it $(9)$ $k_\pi^G(A)$ is a $\pi$-number.}
\smallskip

By Lemma  \ref{ind=inv}

\smallskip
{\it $(10)$ $HA$ leaves invariant each $A$-class of $G$-in\-du\-ced 
$\pi$-Hall subgroups.}
\smallskip

Since  $G\in C_\pi$,

\smallskip
{\it $(11)$ $G$ acts transitively on the set of $A$-classes of
$G$-in\-du\-ced $\pi$-Hall subgroups.}

By $(10)$, the subgroup $HA$ is contained in the kernel of this
action. Now, by $(11)$

\smallskip
{\it $(12)$ $k_\pi^G(A)$ is a $\pi'$-number.}
\smallskip

From $(9)$ and $(12)$ it follows that

\smallskip
{\it $(13)$ $k_\pi^G(A)=1$.}
\smallskip

Now by Lemma \ref{crit}

\smallskip
{\it $(14)$ $HA\in C_\pi$, a contradiction.}
\smallskip

Thus the theorem is proved.

\end{document}